\documentclass[pdflatex,sn-mathphys-ay]{sn-jnl}

\usepackage{graphicx}%
\usepackage{multirow}%
\usepackage{amsmath,amssymb,amsfonts}%
\usepackage{amsthm}%
\usepackage{mathrsfs}%
\usepackage{xcolor}%
\usepackage{textcomp}%
\usepackage{manyfoot}%
\usepackage{booktabs}%
\usepackage{algorithm}%
\usepackage{algorithmicx}%
\usepackage{algpseudocode}%
\usepackage{mathtools}%
\usepackage{array}%
\usepackage{makecell}%

\theoremstyle{thmstyleone}%
\newtheorem{theorem}{Theorem}[section]

\newtheorem{lemma}{Lemma}[section]

\theoremstyle{thmstyletwo}%

\newtheorem{remark}{Remark}[section]

\theoremstyle{thmstylethree}%
\newtheorem{definition}{Definition}[section]

\numberwithin{equation}{section}
\begin{document}

\title[Longest run for $q$-binary trials]{Longest run for $q$-binary trials}

\author*[1,2]{\fnm{Jungtaek} \sur{Oh}}\email{jungtaekoh0191@gmail.com}
\author[2]{\fnm{Nurzhan} \sur{Sapargali}}

\affil*[1]{\orgdiv{Basic Science Research Institute}, \orgname{Pukyong National University}, \orgaddress{\city{Busan}, \postcode{48513}, \country{Republic of Korea}}}

\affil[2]{\orgdiv{Department of Statistics}, \orgname{Ludwig-Maximilians-Universit{\"a}t M{\"u}nchen}, \orgaddress{\street{Ludwigstra{\ss}e 33}, \city{Munich}, \postcode{80539}, \country{Germany}}}

\abstract{We investigate the longest-run distribution for binary trials in which the success probability varies geometrically. Exact formulae are derived for the probability mass function and the cumulative distribution function of the length of the longest success run in a sequence of $q$-binary trials, using enumerative combinatorics. The setting covers independent but not necessarily identically distributed zero--one trials whose success probabilities follow a geometric rule, and the classical independent and identically distributed Bernoulli longest-run distribution is recovered in the limit $q\to 1$. A simulation study confirms accurate maximum likelihood estimation of the base success probability and the geometric parameter $q$, a well-calibrated and powerful boundary-corrected likelihood ratio test of the hypothesis that $q=1$, and close agreement between the analytical cumulative distribution function and its empirical counterpart.}



\maketitle
\section{Introduction}\label{sec:intro}
\begin{figure}[ht]
\centering
\includegraphics[width=0.95\textwidth]{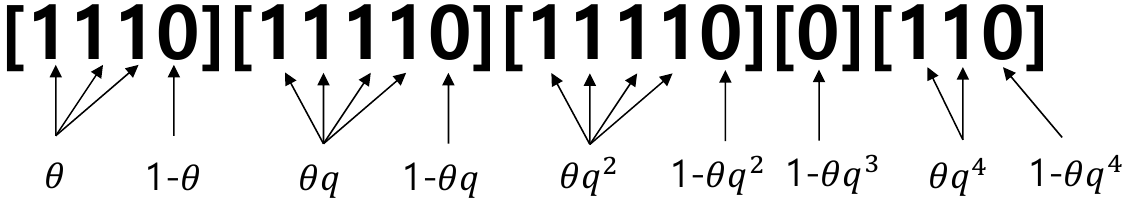}
\caption{A $q$-binary sequence of length $n=18$, partitioned into geometric subsequences with failure probabilities $q_i=1-\theta q^{i-1}$}
\label{fig:qseq}
\end{figure}

\citet{Charalambides2010b} investigated discrete $q$-distributions arising from Bernoulli trials in which the probability of success varies geometrically. Let us consider a sequence $X_{1}$,...,$X_{n}$ consisting of Bernoulli trials that result in either zero (failure) or one (success). The trials within each subsequence, starting immediately after the $(i-1)$st zero and continuing up to the $i$th zero, are independent and share the same probability of failure. For $i>0$, the subsequence of trials that begins after the zero at position $(i-1)$ forms the $i$th geometric sequence. For the $i$th geometric sequence, the probability $q_i$ of zero remains unchanged throughout all trials and is defined as
\noindent
\begin{equation}\label{failureprobofq}
\begin{split}
q_{i}=1-\theta q^{i-1},\quad 0<\theta\leq1,\quad 0<q\leq1.
\end{split}
\end{equation}
\noindent
According to the standard definitions, it follows from equation \eqref{failureprobofq} that, for $0<q<1$, the success probability $1-q_i$ decreases geometrically. Let $S_{j}^{(0)}=\sum_{m=1}^{j}(1-X_{m})$ denote the number of zeros in the first $j$ trials. The probability of observing zeros in the $i$th geometric sequence of trials corresponds to the conditional probability $q_{j,i}$, which represents the conditional probability of occurrence of a zero at trial $j$ given that exactly $i-1$ zeros have occurred in the preceding trials. This relationship can be formally expressed as follows.
\noindent
\begin{equation}\label{failureprobofq2}
\begin{split}
q_{j,i}=P\Big(X_{j}=0\ \bigm|\ S_{j-1}^{(0)}=i-1\Big)=1-\theta q^{i-1},\quad \text{for}\ i\leq j.
\end{split}
\end{equation}
\noindent
To clarify the previous discussion, consider the sequence $111011110111100110$ for $n=18$. The geometric rule indicates that each bracketed subsequence has its own success and failure probabilities, as illustrated in Fig.~\ref{fig:qseq}.

This stochastic model \eqref{failureprobofq} or \eqref{failureprobofq2} has interesting applications, studied as a reliability growth model by \cite{Dubman&Sherman1969}, and applies to a $q$-boson theory in physics by \cite{Jing&Fan1994} and \cite{Jing1994}. The stochastic model \eqref{failureprobofq} is applicable to start-up demonstration tests, functioning as a sequential-intervention model as introduced by \cite{BBV1995}.\\

The stochastic model \eqref{failureprobofq} represents the $q$-analog of the traditional binomial distribution, where the probabilities of zeros change geometrically. This establishes a stochastic framework for independent and identically distributed (IID) trials characterized by a failure probability function  $\pi$  defined as follows:
\noindent
\begin{equation}\label{bernoullifailureprob}
\pi_{j}=P\left(X_{j}=0\right)=1-\theta,\ j>0,
\end{equation}
\noindent
where $0<\theta\leq 1$.

When $q$ approaches $1$, the stochastic model described in \eqref{failureprobofq} simplifies to the IID (Bernoulli) model given in \eqref{bernoullifailureprob}. This occurs because $q_{i}$ converges to $\pi_{i}$ for $i=1,2,\ldots$, or $q_{j,i}$ converges to $1-\theta$ for $i=1,2,\ldots$ and $j=1,2,\ldots.$\\

Numerous scholars have explored discrete $q$-distributions, which arise from the stochastic framework of independent Bernoulli trials. Those interested in a detailed examination and extensive bibliography on this subject should refer to the papers by \citet{charalambides2010a,Charalambides2010b} and the monograph by \citet{Charalambides2016} for further reading.\\

From a Mathematical and Statistical point of view, \cite{Charalambides2016} mentioned in the preface of his book \textit{``It should be noticed that a stochastic model of a sequence of independent Bernoulli trials, in which the probability of success at a trial is assumed to vary with the number of trials and/or the number of successes, is advantageous in the sense that it permits incorporating the experience gained from previous trials and/or successes. If the probability of success at a trial is a very general function of the number of trials and/or the number successes, very little can be inferred from it about the distributions of the various random variables that may be defined on this model. The assumption that the probability of success (or failure) at a trial varies geometrically, with rate (proportion) $q$, leads to the introduction of discrete $q$-distributions''}.\\

In recent decades, the distribution theory of runs and patterns has been extensively developed, inspired by the foundational studies of \citet{philippou1982waiting}, \citet{philippou1983generalized}, and \citet{aki1984discrete}. This area of research has gained attention due to its theoretical significance and its wide-ranging applications in fields such as hypothesis testing, system reliability, quality control, physics, psychology, radar astronomy, molecular biology, computer science, insurance, and finance. These developments have been reviewed by \citet{balakrishnan2003runs} and \citet{fu2003distribution}, along with the references cited therein. Additionally, more recent contributions to this field include the works of \citet{arapis2018}, \citet{eryilmaz2018}, \citet{kong2019}, \citet{makri2019}, \citet{aki2019}, \citet{gera2020acceptance}, and \citet{gera2021runs}.\\

The random variable of interest is $L_{n}^{(1)}$, the length of the longest success run in a sequence of $n$ trials; we write $L_{n}:=L_{n}^{(1)}$ when no confusion arises. Longest runs are important as they represent extreme values within the spectrum of run lengths. Let $R_{n}^{(j)}$ denote the number of runs of type $j$ $(j=0,1)$ and $x_{i}(y_{i})$ represent the length of the $i$th run of successes (failures) in $X_{1}$, $\ldots$,$X_{n}$. In addition, let $L_{n}^{(j)}$ denote the longest run of type $j$($j =0, 1$). To further elaborate, consider a fixed number $n$.
\begin{equation*}
\begin{split}
L_{n}^{(j)}&=\text{max}\Big\{b\ \bigm|\ X_{a+1}=\cdots=X_{a+b}=j\ \wedge\ 0\leq a\leq n-b\Big\},
\quad j=0,1,\\
L_{n}^{(1)}&=\underset{1\leq i \leq R_{n}^{(1)}}{\text{max}}x_{i},\qquad
L_{n}^{(0)}=\underset{1\leq i \leq R_{n}^{(0)}}{\text{max}}y_{i}.
\end{split}
\end{equation*}
For clarity, using an example, we mention that a sequence of ten trials is $0011110110$. We have $R_{10}^{(0)}=3$, $R_{10}^{(1)}=2$, $L_{10}^{(0)}=2$, and $L_{10}^{(1)}=4$.\\
\citet{burr1961longest} studied the exact expression for the distribution of the longest run for IID Bernoulli trials. \citet{fu1994distribution} obtained the distribution of the longest run using a Markov-chain embedding technique. In contrast, \citet{fu1996distribution} developed a finite Markov-chain embedding technique to derive the exact distribution of $L_n$ in sequences of multistate trials. \citet{wendy1996runs} investigated the joint distribution of the length of the longest success run and the number of successes as well as the conditional distribution of the longest-success-run length given the number of successes in sequences of Bernoulli trials. Asymptotic results for the longest-run statistics in a multistate Markov chain were explored by \citet{vaggelatou2003length}. Meanwhile, \citet{makri2007shortest} derived the exact probability function of the longest-run statistics for linear and circular sequences through combinatorial analysis. Numerous researchers have contributed to studying longest-run statistics in various scenarios, including Bernoulli trials, Markov-dependent trials, and exchangeable binary trials. For further insights, refer to works such as those by \citet{philippou1985longest,philippou1986successes,philippou1990longest},  \citet{antzoulakos1997probability,antzoulakos1998longest}, \citet{charalambides1994success}, \citet{eryilmaz2005distribution,eryilmaz2005longest,eryilmaz2006some}, and \citet{muselli1996simple}.\\
As mentioned by \cite{johnson2005univariate}, ``\textit{The formulas for the probabilities for these distributions are not at first sight very illuminating; also in some cases, there is more than one correct published expression for the probabilities}.'' \citet{fu1994distribution} mentioned that some formulas are complicated and hard to compute. In this paper, we derive explicit probability functions for the longest success run in binary trials whose success probabilities vary according to a geometric rule. Our results provide a direct longest-run formulation based on enumerative combinatorics and are complemented by numerical verification and statistical inference for the model parameters.\\
Herein, we study the distribution of the length of the longest success run, with the probability of ones varying according to a geometric rule. The paper is organized as follows. In Section 2 we introduce basic definitions and necessary notations that will be useful throughout this article.
In Section \ref{sec:longes_run_qbinary} we shall study the probability mass and cumulative distribution functions of the length of the longest success run. 
We derive the exact probability function and the cumulative distribution function (CDF) of the longest-success-run length through combinatorial analysis, and illustrate these results numerically for indicative parameter values.
In Section \ref{sec:sim}, we conduct a simulation study to assess the accuracy of the maximum likelihood estimator of $(\theta, q)$, to evaluate the size and power of a boundary-corrected likelihood ratio test of $H_0: q=1$, and to compare the derived analytical cumulative distribution function of the longest-success-run length with its empirical counterpart.


\section{Terminology and notation}\label{sec:notation}
We first recall some definitions, notations, and known results, which will be used in this paper. Throughout the paper, we suppose that $0<q\leq 1$. First, we introduce the following notation.
\begin{itemize}
\item $L_n^{(1)}:$ the length of the longest run of successes in $X_{1}, X_{2}, \ldots, X_{n}$;
\item $L_n^{(0)}:$ the length of the longest run of failures in $X_{1}, X_{2}, \ldots, X_{n}$;
\item $S_{n}:$ the total number of successes in $X_{1}, X_{2}, \ldots, X_{n}$;
\item $F_{n}:$ the total number of failures in $X_{1}, X_{2}, \ldots, X_{n}$.
\end{itemize}
\noindent
We write $(a;q)_{0}=1$ and $(a;q)_{n}=\prod_{k=0}^{n-1}(1-aq^{k})$ for $n\geq 1$. For $0<q<1$ and integers $n\geq r\geq 0$, the $q$-binomial coefficient (Gaussian polynomial) is
\begin{equation*}
\begin{bmatrix}
n\\
r
\end{bmatrix}_{q}
=\frac{(q;q)_{n}}{(q;q)_{r}(q;q)_{n-r}},
\end{equation*}
and the case $q=1$ is understood in the limiting sense $\lim_{q\to 1}\begin{bmatrix}n\\r\end{bmatrix}_{q}=\binom{n}{r}$.
\noindent
Let us consider again a sequence of independent geometric sequences of trials with probability of failure at the $i$th geometric sequence of trials given by \eqref{failureprobofq} or \eqref{failureprobofq2}. We now focus on the number of successes in a given number of trials in this stochastic model.
\begin{definition}
Let $Z_{n}$ be the number of successes in a sequence of $n$ independent Bernoulli trials, with probability of failure at the $i$th geometric sequence of trials given by \eqref{failureprobofq} or \eqref{failureprobofq2}. The distribution of the random variable $Z_{n}$ is called $q$-binomial distribution, with parameters $n$, $\theta$, and $q$.
\end{definition}
\noindent
We now introduce a $q$-analog of the binomial distribution. The probability function of the number $Z_{n}$ of successes in $n$ trials $X_{1}, \ldots, X_{n}$ is given by
\begin{equation}\label{q-binomialdist}
\begin{split}
P_{q,\theta}\{Z_{n}=r\}=
\begin{bmatrix}
n\\
r
\end{bmatrix}_{q}
\theta^{r}\prod_{i=1}^{n-r}(1-\theta q^{i-1}),
\end{split}
\end{equation}\\
for $r = 0, 1, \ldots, n$ and $0<q\leq 1$. The distribution is called a $q$-binomial distribution. For $q \rightarrow 1$, because
\begin{equation*}
\begin{split}
\lim\limits_{q\rightarrow 1}
\begin{bmatrix}
n\\
r
\end{bmatrix}_{q}
={n \choose r}
\end{split}
\end{equation*}
so that the $q$-binomial distribution converges to the usual binomial distribution as $q \rightarrow 1$, as follows
\begin{equation}
P_{\theta}\left(Z_{n}=r\right)={n\choose r}\theta^{r}(1-\theta)^{n-r},\ r=0,1,\ldots,n,
\end{equation}
with parameters $n$ and $\theta$.
The $q$-binomial distribution studied by \citet{charalambides2010a,Charalambides2016} is connected with the $q$-Bernstein polynomials. \citet{Jing1994} introduced the probability function \eqref{q-binomialdist} as a $q$-deformed binomial distribution and also derived a recurrence relation for its probability distribution.
\noindent
In the sequel, $P_{q,\theta}(\cdot)$ and $P_{\theta}(\cdot)$ denote probabilities related with the stochastic model \eqref{failureprobofq} and \eqref{bernoullifailureprob}, respectively.

\section{Longest run for $q$-binary trials}
\label{sec:longes_run_qbinary}

In this section, we delve into the distribution of the length of the longest success run for binary (zero and one) trials, where the probability of encountering ones follows a geometric rule. First, we derive the exact expression for the PMF $P_{q,\theta}\left(L_{n}^{(1)}= k\right)$, representing the distribution of the length $L_{n}^{(1)}$ of the longest success run. Subsequently, we obtain the cumulative distribution function $P_{q,\theta}\left(L_{n}^{(1)}\leq k\right)$ for $L_{n}^{(1)}$, the longest-success-run length. These results are derived through enumerative combinatorics.

\subsection{Closed-form expression for the PMF of $L_{n}^{(1)}$}

\noindent
We now present a definition and a lemma that will be used in the proofs of the theorems in the sequel.
\noindent
\begin{definition}
For $0<q\leq1$, we define

\begin{equation*}
\begin{split}
U_{q}^{k}(r,s,t)={\sum_{x_{1},\ldots,x_{r}}} q^{x_{2}+2x_{3}+\cdots+(r-1)x_{r}},\
\end{split}
\end{equation*}
\noindent
where the summation is over all integers $x_1,\ldots,x_{r}$ satisfying
\noindent

\begin{equation*}
\begin{split}
0\leq x_{j}\leq k\ \text{for}\ j=1,\ldots,r,
\end{split}
\end{equation*}
\noindent
\begin{equation*}
\begin{split}
x_{1}+\cdots+x_{r}=s,\ \text{and}
\end{split}
\end{equation*}
\noindent
\begin{equation*}
\begin{split}
\delta_{x_{1},k}+\cdots+\delta_{x_{r},k}=t.
\end{split}
\end{equation*}

\end{definition}

\noindent
The following gives a recurrence relation useful for the computation of $U_{q}^{k}(r,s,t)$.\\
\noindent

\begin{lemma}
\label{lemma3.1-n}
For $0<q\leq1$, $U_{q}^{k}(r,s,t)$ obeys the following recurrence relation.
\noindent
\begin{equation*}
\begin{split}
U_{q}&^{k}(r,s,t)\\
&=\left\{
  \begin{array}{ll}
    \sum_{a=0}^{k-1}q^{a(r-1)}U_{q}^{k}(r-1,s-a,t)+\\
    \quad \quad\ q^{k(r-1)}U_{q}^{k}(r-1,s-k,t-1) & \text{for}\ r>1,\ tk\leq s\leq rk\ \text{and}\ t\leq r, \\
    1 & \text{for}\ r=1,\ s=k,\ \text{and}\ t=1,\\
    & \text{or}\ r=1,\ 0\leq s<k,\ \text{and}\ t=0,\\
    0 & \text{otherwise.}\\
  \end{array}
\right.
\end{split}
\end{equation*}
\end{lemma}
\begin{proof}
For $r > 1$, $tk\leq s\leq rk$ and $t\leq r$, we observe that since $x_{r}$ can assume the values $0,\ldots,k$, then $U_{q}^{k}(r,s,t)$ can be written as follows.
\noindent
\begin{equation*}
\begin{split}
U_{q}^{k}(r,s,t)=&\sum_{x_{r}=0}^{k-1}\ q^{x_{r}(r-1)} {\hspace{0.5cm}\sum_{\substack{x_{1}+\cdots+x_{r-1}=s-x_{r}\\ x_{1},\ldots,x_{r-1} \in \{0,1,\ldots,k\}\\\delta_{x_{1},k}+\cdots+\delta_{x_{r-1},k}=t}}}\quad
 q^{x_{2}+2x_{3}+\cdots+(r-2)x_{r-1}}\\
 &+q^{k(r-1)}{\hspace{0.5cm}\sum_{\substack{x_{1}+\cdots+x_{r-1}=s-k\\ x_{1},\ldots,x_{r-1} \in \{0,1,\ldots,k\}\\\delta_{x_{1},k}+\cdots+\delta_{x_{r-1},k}=t-1}}}\quad
\ q^{x_{2}+2x_{3}+\cdots+(r-2)x_{r-1}}\\
\end{split}
\end{equation*}
\noindent
Using simple algebraic arguments to simplify, then $U_{q}^{k}(r,s,t)$ can be rewritten as follows.
\noindent
\begin{equation*}
\begin{split}
U_{q}^{k}(r,s,t)=&\sum_{a=0}^{k-1}q^{a(r-1)}U_{q}^{k}(r-1,s-a,t)+ q^{k(r-1)}U_{q}^{k}(r-1,s-k,t-1).
\end{split}
\end{equation*}
\noindent
The other cases are obvious and thus the proof is completed.
\end{proof}

\begin{remark}
{\rm
We observe that $U_1^{k}(r,s,t)$ is the number of integer solutions $(x_{1},\ \ldots,\ x_{r})$ of
\noindent

\begin{equation*}
\begin{split}
0\leq x_{j}\leq k\ \text{for}\ j=1,\ldots,r,
\end{split}
\end{equation*}
\noindent
\begin{equation*}
\begin{split}
x_{1}+\cdots+x_{r}=s,\ \text{and}
\end{split}
\end{equation*}
\noindent
\begin{equation*}
\begin{split}
\delta_{x_{1},k}+\cdots+\delta_{x_{r},k}=t.
\end{split}
\end{equation*}
\noindent
which is
\noindent
\begin{equation*}
\begin{split}
U_1^{k}(r,s,t)={r \choose t}C(s-tk,\ r-t,\ k-1),
\end{split}
\end{equation*}
\noindent
denotes the number of ways of distributing $s$ identical balls into $r$ different cells while ensuring that each of $t$ cells receives exactly $k$ balls, with the remaining $s-tk$ balls placed in the remaining $r-t$ cells such that no cell receives more than $k-1$ balls. In addition, $C(a,\ b,\ c)$ represents the total number of integer solutions to the equation $x_{1}+x_{2}+\cdots+x_{b}=a$, subject to the condition $0\leq x_{i}\leq c$ for $i=1,2,\ldots,b$. Alternatively, $C(a, b, c)$ indicates the number of ways of distributing $a$ identical balls into $b$ different cells while ensuring that each cell is occupied by at most $c$ balls. This number can be expressed as follows: 
\noindent
\begin{equation*}
\begin{split}
C(a,\ b,\ c)=\sum_{j=0}^{b}(-1)^{j}{b \choose j}{a-(c+1)j+b-1 \choose b-1}.
\end{split}
\end{equation*}
\noindent
Here $C(a,b,c)$ equals the displayed sum if $a>0$, equals $1$ if $a=0$, and equals $0$ otherwise (see, e.g., \citet{charalambides2002enumerative}).
}
\end{remark}
\noindent
We shall study the distribution of the length $L_{n}^{(1)}$ of the longest success run. The probability mass function $P_{q,\theta}\Big(L_{n}^{(1)}= k\Big)$ is obtained by the following theorem. It is evident that
\noindent
\begin{equation*}
P_{q,\theta}\left(L_{n}^{(1)}=0\right)={\prod_{j=1}^{n}}\left(1-\theta q^{j-1}\right)\ \text{for}\ k=0
\end{equation*}
\noindent
and so we shall focus on determining the probability mass function for $1\leq k \leq n$.
\noindent
\begin{theorem}
\label{thm:3.1-n}
The PMF $P_{q,\theta}\Big(L_{n}^{(1)}= k\Big)$, for $k=1,\ldots,n$, is given by
\begin{equation*}
\begin{split}
P_{q,\theta}\Big(L_{n}^{(1)}= k\Big)=&\sum_{y=[n/(k+1)]}^{n-k}\theta^{n-y}\ {\prod_{j=1}^{y}}\ \left(1-\theta q^{j-1}\right)\ \sum_{i=1}^{y+1} U_{q}^{k}(y+1,\ n-y,\ i).
\end{split}
\end{equation*}
\end{theorem}
\begin{proof}
We partition the event $\left\{L_{n}^{(1)}= k\right\}$ into disjoint events given by $F_{n}=y,$ for $y=[n/(k+1)], \ldots,n-k$. Adding the probabilities, we have
\noindent
\begin{equation*}
\begin{split}
P_{q,\theta}\Big(L_{n}^{(1)}= k\Big)=\sum_{y=[n/(k+1)]}^{n-k}P_{q,\theta}\Big(L_{n}^{(1)}= k\ \wedge\ F_{n}=y\Big).\\
\end{split}
\end{equation*}
\noindent
We focus on the event $\left\{L_{n}^{(1)}= k\ \wedge\ F_{n}=y\right\}$. A typical element of this event, denoted as $\Big\{L_{n}^{(1)}= k\ \wedge\ F_{n}=y\Big\}$, consists of an ordered sequence comprising $n-y$ successes and $y$ failures such that the length of the longest success run is equal to $k$. The number of such sequences can be derived as follows. First, we distribute the $y$ failures, forming $y+1$ cells. Next, we distribute the $n-y$ successes among these $y+1$ distinguishable cells in a way that satisfies the following conditions: among these cells, denoted as $i$ ($i=1,\ldots,\text{min}\{y+1,[(n-y)/k]\})$, exactly $k$ successes occupy each of the specified cells. The remaining $n-y-ik$ successes can then be distributed among the remaining $y+1-i$ cells, ensuring no cell receives more than $k-1$ successes.
Thus, the event $\left\{L_{n}^{(1)}= k\ \wedge\ F_{n}=y\right\}$ is defined as follows.
\noindent
\begin{equation*}
\begin{split}
\underbrace{1\ldots 1}_{0\leq x_{1}\leq k}0\underbrace{1\ldots 1}_{0\leq x_{2}\leq k}0\underbrace{1\ldots 1}_{0\leq x_{3}\leq k}0 \ldots 0\underbrace{1\ldots 1}_{0\leq x_{y}\leq k}0\underbrace{1\ldots 1}_{0\leq x_{y+1}\leq k},\\
\end{split}
\end{equation*}
\noindent
with $y$ $0$'s and $n-y$ $1$'s, where $x_{j}$ $(j=1,\ldots,y+1)$ represents a length of a run of $1$'s. And all integers $x_{1},\ldots,x_{y+1}$ satisfy the conditions
\noindent
\begin{equation*}
\begin{split}
0 \leq x_j \leq k \mbox{\ for\ } j=1,...,y+1,\  x_1+\cdots +x_{y+1} = n-y,\mbox{\ and}
\end{split}
\end{equation*}
\noindent
\begin{equation*}
\begin{split}
\delta_{x_1,k}+\delta_{x_2,k}+\cdots+\delta_{x_{y+1},k}=i.
\end{split}
\end{equation*}
\noindent
Then the probability of the event $\left\{L_{n}^{(1)}= k\ \wedge\ F_{n}=y\right\}$ is given by
\noindent
\begin{equation*}
\begin{split}
P_{q,\theta}\Big(L_{n}^{(1)}= k\ \wedge\ F_{n}=y\Big)&\\
=\sum_{i=1}^{y+1}\Bigg[\bigg\{\sum_{\substack{x_1+\cdots +x_{y+1} = n-y\\ x_{1},\ldots,x_{y+1} \in \{0,1,\ldots,k\}\\\delta_{x_1,k}+\cdots+\delta_{x_{y+1},k}=i}}&\Big(\theta q^{0}\Big)^{x_{1}}\Big(1-\theta q^{0}\Big)\Big(\theta q^{1}\Big)^{x_{2}}\Big(1-\theta q^{1}\Big)\times\\
&\quad \quad \quad \quad\quad \quad\quad \quad \quad  \vdots\\
& \Big(1-\theta q^{y-2}\Big)\Big(\theta q^{y-1}\Big)^{x_{y}} \Big(1-\theta q^{y-1}\Big)\Big(\theta q^{y}\Big)^{x_{y+1}}\bigg\}\Bigg].\\
\end{split}
\end{equation*}
\noindent
Using simple exponentiation algebra arguments to simplify,
\noindent
\begin{equation*}
\begin{split}
P_{q,\theta}\Big(&L_{n}^{(1)}= k\ \wedge\ F_{n}=y\Big)=\\
&\theta^{n-y}\ {\prod_{j=1}^{y}}\ \left(1-\theta q^{j-1}\right)\ \sum_{i=1}^{y+1}\ \sum_{\substack{x_1+\cdots +x_{y+1} = n-y\\ x_{1},\ldots,x_{y+1} \in \{0,1,\ldots,k\}\\\delta_{x_1,k}+\cdots+\delta_{x_{y+1},k}=i}}q^{x_{2}+2x_{3}+\cdots+(y-1)x_{y}+yx_{y+1}}.\\
\end{split}
\end{equation*}
\noindent
Using Lemma \ref{lemma3.1-n}, we can rewrite as follows.
\noindent
\begin{equation*}
\begin{split}
P_{q,\theta}\Big(L_{n}^{(1)}= k\ \wedge\ F_{n}=y\Big)=\theta^{n-y}\ {\prod_{j=1}^{y}}\ \left(1-\theta q^{j-1}\right)\ \sum_{i=1}^{y+1} U_{q}^{k}(y+1,\ n-y,\ i),\\
\end{split}
\end{equation*}\\
\noindent
where
\noindent
\begin{equation*}
\begin{split}
U_{q}^{k}(y+1,\ n-y,\ i)=\sum_{\substack{x_1+\cdots +x_{y+1} = n-y\\ x_{1},\ldots,x_{y+1} \in \{0,1,\ldots,k\}\\\delta_{x_1,k}+\cdots+\delta_{x_{y+1},k}=i}}q^{x_{2}+2x_{3}+\cdots+(y-1)x_{y}+yx_{y+1}}.
\end{split}
\end{equation*}
\noindent
Therefore we can compute the probability of the event $\left\{L_{n}^{(1)}= k\right\}$ as follows.
\noindent
\begin{equation*}
\begin{split}
P_{q,\theta}\Big(L_{n}^{(1)}= k\Big)=&\sum_{y=[n/(k+1)]}^{n-k}\theta^{n-y}\ {\prod_{j=1}^{y}}\ \left(1-\theta q^{j-1}\right)\ \sum_{i=1}^{y+1} U_{q}^{k}(y+1,\ n-y,\ i).\\
\end{split}
\end{equation*}
\noindent
Thus, proof is completed.
\end{proof}
\noindent
It is worth mentioning here that the PMF $P_{q,\theta}\Big(L_{n}^{(1)}= k\Big)$ approaches the PMF $P_{\theta}\Big(L_{n}^{(1)}= k\Big)$ of the IID model as $q$ tends to $1$. The details are presented in the following remark.
\noindent
\begin{remark}
{\rm For $q=1$, the PMF $P_{q,\theta}\Big(L_{n}^{(1)}= k\Big)$ reduces to the PMF $P_{\theta}\Big(L_{n}^{(1)}= k\Big)$. Then $P_{\theta}\Big(L_{n}^{(1)}= 0\Big)=(1-\theta)^{n}$, and for $1\leq k \leq n$ the PMF is given by
\noindent
\begin{equation*}
\begin{split}
P_{\theta}\Big(L_{n}^{(1)}= k\Big)=&\sum_{y=[n/(k+1)]}^{n-k}\theta^{n-y}(1-\theta)^{y}\ \sum_{i=1}^{y+1} U_{1}^{k}(y+1,\ n-y,\ i)\\
=&\sum_{y=[n/(k+1)]}^{n-k}\theta^{n-y}(1-\theta)^{y}\ \sum_{i=1}^{y+1} {y+1 \choose i}C(n-y-ik,\ y+1-i,\ k-1),
\end{split}
\end{equation*}
\noindent
where
\noindent
\begin{equation*}
\begin{split}
C(a,\ b,\ c)=\sum_{j=0}^{b}(-1)^{j}{b \choose j}{a-(c+1)j+b-1 \choose b-1}.
\end{split}
\end{equation*}
\noindent
Here $C(a,b,c)$ equals the displayed sum if $a>0$, equals $1$ if $a=0$, and equals $0$ otherwise (see, e.g., \citet{charalambides2002enumerative}).
}
\end{remark}

\subsection{Closed-form expression for the CDF of $L_{n}^{(1)}$}
\noindent
We now present a definition and a lemma that will be used in the proofs of the theorems in the sequel.
\noindent
\begin{definition}
For $0<q\leq1$, we define
\noindent
\begin{equation*}
\begin{split}
V_{q}(r,s,k)={\sum_{x_{1},\ldots,x_{r}}} q^{x_{2}+2x_{3}+\cdots+(r-1)x_{r}},\
\end{split}
\end{equation*}
\noindent
where the summation is over all integers $x_1,\ldots,x_{r}$ satisfying
\noindent
\begin{equation*}
\begin{split}
0\leq x_{j}\leq k\quad \text{for}\quad j=1,\ldots,r,\ \text{and}
\end{split}
\end{equation*}
\noindent
\begin{equation*}
\begin{split}
x_{1}+\cdots+x_{r}=s
\end{split}
\end{equation*}
\noindent
\end{definition}
\noindent
The following gives a recurrence relation useful for the computation of $V_{q}(r,s,k)$.\\
\noindent
\begin{lemma}
\label{lemma3.2-n}
For $0<q\leq1$, $V_{q}(r,s,k)$ obeys the following recurrence relation.
\noindent
\begin{equation*}
\begin{split}
V_{q}&(r,s,k)=\left\{
  \begin{array}{ll}
    \sum_{a=0}^{k}q^{a(r-1)}V_{q}(r-1,s-a,k) & \text{for}\ r>1,\ \text{and}\ 0\leq s\leq kr\\
     1 & \text{for}\ r=1,\ \text{and}\ 0\leq s\leq k,\ \text{or}\\
    0 & \text{otherwise.}\\
  \end{array}
\right.
\end{split}
\end{equation*}
\end{lemma}
\begin{proof}
For $r > 1$ and $0\leq s\leq kr$, we observe that since $x_{r}$ can assume the values $0,\ldots,k$, then $V_{q}(r,s,k)$ can be written as follows.
\noindent
\begin{equation*}
\begin{split}
V_{q}(r,s,k)=&\sum_{x_{r}=0}^{k}\ q^{x_{r}(r-1)} {\hspace{0.5cm}\sum_{\substack{x_{1}+\cdots+x_{r-1}=s-x_{r}\\ x_{1},\ldots,x_{r-1} \in \{0,1,\ldots,k\}}}}\quad
 q^{x_{2}+2x_{3}+\cdots+(r-2)x_{r-1}}
\end{split}
\end{equation*}
\noindent
Using simple algebraic arguments to simplify, then $V_{q}(r,s,k)$ can be rewritten as follows.
\noindent
\begin{equation*}
\begin{split}
V_{q}&(r,s,k)=\sum_{a=0}^{k}q^{a(r-1)}V_{q}(r-1,s-a,k).
\end{split}
\end{equation*}
\noindent
The other cases are obvious and thus the proof is completed.
\end{proof}

\begin{remark}
{\rm
We observe that $V_1(r,s,k)$ is the number of integer solutions $(x_{1},\ \ldots,\ x_{r})$ of
\noindent
\begin{equation*}
\begin{split}
0\leq x_{j}\leq k\quad \text{for}\quad j=1,\ldots,r,\ \text{and}
\end{split}
\end{equation*}
\noindent
\begin{equation*}
\begin{split}
x_{1}+\cdots+x_{r}=s
\end{split}
\end{equation*}
\noindent
which is
\noindent
\begin{equation*}
\begin{split}
V_1(r,s,k)=C(s,\ r,\ k),
\end{split}
\end{equation*}
\noindent
where $C(a,\ b,\ c)$ denotes the total number of integer solutions $x_{1}+x_{2}+\cdots+x_{b}=a$ such that $0\leq x_{i}\leq c$ for $i=1,2,\ldots,b$. Alternatively, it is the number of ways of distributing $a$ identical balls into $b$ different cells so that each cell is occupied by at most $c$ balls. The number can be expressed as
\noindent
\begin{equation*}
\begin{split}
C(a,\ b,\ c)=\sum_{j=0}^{b}(-1)^{j}{b \choose j}{a-(c+1)j+b-1 \choose b-1}.
\end{split}
\end{equation*}
\noindent
Here $C(a,b,c)$ equals the displayed sum if $a>0$, equals $1$ if $a=0$, and equals $0$ otherwise (see, e.g., \citet{charalambides2002enumerative}).
}
\end{remark}
\noindent
We shall study the cumulative distribution function $P_{q,\theta}\Big(L_{n}^{(1)}\leq k\Big)$ of the length $L_{n}^{(1)}$ of the longest success run. The CDF $P_{q,\theta}\Big(L_{n}^{(1)}\leq k\Big)$ is obtained by the following theorem.
\noindent
\begin{theorem}
\label{thm:3.2-n}
The cumulative distribution function $P_{q,\theta}\Big(L_{n}^{(1)}\leq k\Big)$ of the length $L_{n}^{(1)}$ of the longest success run is given by
\begin{equation*}
\begin{split}
P_{q,\theta}\Big(L_{n}^{(1)}\leq k\Big)=&\sum_{y=[n/(k+1)]}^{n}\theta^{n-y}\ {\prod_{j=1}^{y}}\ \left(1-\theta q^{j-1}\right)\  V_{q}(y+1,\ n-y,\ k).
\end{split}
\end{equation*}
\end{theorem}

\begin{proof}
We partition the event $\left\{L_{n}^{(1)}\leq k\right\}$ into disjoint events given by $F_{n}=y,$ for $y=[n/(k+1)], \ldots,n$. Adding the probabilities we have
\noindent
\begin{equation*}
\begin{split}
P_{q,\theta}\Big(L_{n}^{(1)}\leq k\Big)=\sum_{y=[n/(k+1)]}^{n}P_{q,\theta}\Big(L_{n}^{(1)}\leq k\ \wedge\ F_{n}=y\Big).\\
\end{split}
\end{equation*}
\noindent
We focus on the event $\left\{L_{n}^{(1)}\leq k\ \wedge\ F_{n}=y\right\}$. A typical element of this event, denoted as $\Big\{L_{n}^{(1)}\leq k\ \wedge\ F_{n}=y\Big\}$, consists of an ordered sequence comprising $n-y$ successes and $y$ failures such that the length of the longest success run is less than or equal to $k$. The number of such sequences can be derived as follows. First, we distribute the $y$ failures, forming $y+1$ cells. Next, we distribute the $n-y$ successes among these $y+1$ distinguishable cells, ensuring that each cell is occupied by at most $k$ successes. Thus, the event $\left\{L_{n}^{(1)}\leq k\ \wedge\ F_{n}=y\right\}$ is defined as follows:
\noindent
\begin{equation*}
\begin{split}
\underbrace{1\ldots 1}_{0\leq x_{1}\leq k}0\underbrace{1\ldots 1}_{0\leq x_{2}\leq k}0\underbrace{1\ldots 1}_{0\leq x_{3}\leq k}0 \ldots 0\underbrace{1\ldots 1}_{0\leq x_{y}\leq k}0\underbrace{1\ldots 1}_{0\leq x_{y+1}\leq k},\\
\end{split}
\end{equation*}
\noindent
with $y$ $0$'s and $n-y$ $1$'s, where $x_{j}$ $(j=1,\ldots,y+1)$ represents a length of a run of $1$'s. In addition, all integers $x_{1},\ldots,x_{y+1}$ satisfy the conditions
\noindent
\begin{equation*}
\begin{split}
0 \leq x_j \leq k \mbox{\ for\ } j=1,...,y+1,\  \mbox{\ and}\ x_1+\cdots +x_{y+1} = n-y,
\end{split}
\end{equation*}
\noindent
Then the probability of the event $\left\{L_{n}^{(1)}\leq k\ \wedge\ F_{n}=y\right\}$ is given by
\noindent
\begin{equation*}
\begin{split}
P_{q,\theta}\Big(L_{n}^{(1)}\leq k\ \wedge\ F_{n}=y\Big)&\\
=\sum_{\substack{x_1+\cdots +x_{y+1} = n-y\\ x_{1},\ldots,x_{y+1} \in \{0,1,\ldots,k\}}}&\Big(\theta q^{0}\Big)^{x_{1}}\Big(1-\theta q^{0}\Big)\Big(\theta q^{1}\Big)^{x_{2}}\Big(1-\theta q^{1}\Big)\times\\
&\quad \quad \quad \quad\quad \quad\quad \quad \quad  \vdots\\
& \Big(1-\theta q^{y-2}\Big)\Big(\theta q^{y-1}\Big)^{x_{y}} \Big(1-\theta q^{y-1}\Big)\Big(\theta q^{y}\Big)^{x_{y+1}}.\\
\end{split}
\end{equation*}
\noindent
Using simple exponentiation algebra arguments to simplify,
\noindent
\begin{equation*}
\begin{split}
P_{q,\theta}\Big(&L_{n}^{(1)}\leq k\ \wedge\ F_{n}=y\Big)=\\
&\theta^{n-y}\ {\prod_{j=1}^{y}}\ \left(1-\theta q^{j-1}\right)\ \sum_{\substack{x_1+\cdots +x_{y+1} = n-y\\ x_{1},\ldots,x_{y+1} \in \{0,1,\ldots,k\}}}q^{x_{2}+2x_{3}+\cdots+(y-1)x_{y}+yx_{y+1}}.\\
\end{split}
\end{equation*}
\noindent
Using Lemma \ref{lemma3.2-n}, we can rewrite as follows.
\noindent
\begin{equation*}
\begin{split}
P_{q,\theta}\Big(L_{n}^{(1)}\leq k\ \wedge\ F_{n}=y\Big)=\theta^{n-y}\ {\prod_{j=1}^{y}}\ \left(1-\theta q^{j-1}\right)\  V_{q}(y+1,\ n-y,\ k),\\
\end{split}
\end{equation*}\\
\noindent
where
\noindent
\begin{equation*}
\begin{split}
V_{q}(y+1,\ n-y,\ k)=\sum_{\substack{x_1+\cdots +x_{y+1} = n-y\\ x_{1},\ldots,x_{y+1} \in \{0,1,\ldots,k\}}}q^{x_{2}+2x_{3}+\cdots+(y-1)x_{y}+yx_{y+1}}.
\end{split}
\end{equation*}
\noindent
Therefore we can compute the probability of the event $\left\{L_{n}^{(1)}\leq k\right\}$ as follows.
\noindent
\begin{equation*}
\begin{split}
P_{q,\theta}\Big(L_{n}^{(1)}\leq k\Big)=&\sum_{y=[n/(k+1)]}^{n}\theta^{n-y}\ {\prod_{j=1}^{y}}\ \left(1-\theta q^{j-1}\right)\ V_{q}(y+1,\ n-y,\ k).\\
\end{split}
\end{equation*}
\noindent
Thus proof is completed.
\end{proof}
It is worth mentioning here that the CDF $P_{q,\theta}\Big(L_{n}^{(1)}\leq k\Big)$ approaches the CDF $P_{\theta}\Big(L_{n}^{(1)}\leq k\Big)$ of the IID model in the limit as $q$ tends to $1$. The details are presented in the following remark.
\begin{remark}
{\rm For $q=1$, the CDF $P_{q,\theta}\Big(L_{n}^{(1)}\leq k\Big)$ reduces to the CDF $P_{\theta}\Big(L_{n}^{(1)}\leq k\Big)$, which is given by
\noindent
\begin{equation*}
\begin{split}
P_{\theta}\Big(L_{n}^{(1)}\leq k\Big)=&\sum_{y=[n/(k+1)]}^{n}\theta^{n-y}(1-\theta)^{y}\ V_{1}(y+1,\ n-y,\ k)\\
=&\sum_{y=[n/(k+1)]}^{n}\theta^{n-y}(1-\theta)^{y}\ C(n-y,\ y+1,\ k),
\end{split}
\end{equation*}
\noindent
where
\noindent
\begin{equation*}
\begin{split}
C(a,\ b,\ c)=\sum_{j=0}^{b}(-1)^{j}{b \choose j}{a-(c+1)j+b-1 \choose b-1}.
\end{split}
\end{equation*}
\noindent
Here $C(a,b,c)$ equals the displayed sum if $a>0$, equals $1$ if $a=0$, and equals $0$ otherwise (see, e.g., \citet{charalambides2002enumerative}).
}
\end{remark}
\begin{remark}
Let $G_{n,k}$ denote the number of success runs of length at least $k$
in a $q$-sequence, as studied by \citet{makri2016runs}. Since
\[
\{L_n^{(1)}\leq k\}=\{G_{n,k+1}=0\},
\]
we have
\[
P_{q,\theta}(L_n^{(1)}\leq k)
=P_{q,\theta}(G_{n,k+1}=0).
\]
Thus, the CDF of the longest success run can also be expressed as a
zero-count probability from the distribution of $G_{n,k+1}$. The
expression in Theorem~\ref{thm:3.2-n} therefore provides an alternative
closed-form representation of this CDF directly in terms of the
longest-run statistic.
\end{remark}

\begin{remark}
Using CDF in Theorem \ref{thm:3.2-n} an alternative expression for the $r$th raw moment (alternatively the $r$th uncorrected moment or the $r$th moment about zero) is given by
\noindent
\begin{equation} \label{eq:moments}
    \mu_{q}(r;n;\theta)=E_{q,\theta}[(L_n^{(1)})^r] =\sum_{j=0}^{n-1}\big[(j+1)^r-j^r\big]\big[1-P_{q,\theta}(L_n^{(1)}\leq j)\big],
\end{equation}
which is directly computable, as $n$ is finite.
\end{remark}

\subsection{Numerics}
\label{subsec:numerics}

Some numerics are next presented, computed by implementation of the PMF and CDF of $L_n^{(1)}$ obtained in Theorems~\ref{thm:3.1-n} and~\ref{thm:3.2-n}.
Table~\ref{tab:pmf} depicts the PMF $P_{q,\theta}(L_n^{(1)}=k)$, the mean $E_{q, \theta}[L_n^{(1)}]$ and variance $\mathrm{Var}_{q, \theta}(L_n^{(1)})$, while Table~\ref{tab:cdf} depicts the CDF $P_{q,\theta}(L_n^{(1)}\leq k)$, for $n\in\{3,5\}$ and for the indicative values $\theta\in\{0.25,0.75\}$, $q\in\{0.25,0.75\}$.

\begin{sidewaystable}
    \centering
    \caption{PMF $P_{q,\theta}(L_n^{(1)}=k)$, mean $E_{q, \theta}[L_n^{(1)}]$ and variance
    $\mathrm{Var}_{q, \theta}(L_n^{(1)})$, for $n=3,5$ and $(\theta,q)\in\{0.25,0.75\}^2$.}
    \label{tab:pmf}
    \begin{tabular}{ccccccccccc}
        \toprule
        \multirow{2}{*}{$n$} & \multirow{2}{*}{$\theta$} & \multirow{2}{*}{$q$} &
        \multicolumn{6}{c}{$P_{q,\theta}(L_n^{(1)}=k)$} &
        \multirow{2}{*}{$E_{q, \theta}[L_n^{(1)}]$} & \multirow{2}{*}{$\mathrm{Var}_{q, \theta}(L_n^{(1)})$} \\
        \cmidrule{4-9}
         & & & $k=0$ & $k=1$ & $k=2$ & $k=3$ & $k=4$ & $k=5$ & & \\
        \midrule
        3 & 0.25 & 0.25 & 0.6921 & 0.2424 & 0.0498 & 0.0156 & --     & --     & 0.3889 & 0.4310 \\
        3 & 0.25 & 0.75 & 0.5237 & 0.3875 & 0.0732 & 0.0156 & --     & --     & 0.5808 & 0.4837 \\
        3 & 0.75 & 0.25 & 0.1936 & 0.2351 & 0.1494 & 0.4219 & --     & --     & 1.7996 & 1.3912 \\
        3 & 0.75 & 0.75 & 0.0632 & 0.2952 & 0.2197 & 0.4219 & --     & --     & 2.0002 & 0.9700 \\
        \addlinespace
        5 & 0.25 & 0.25 & 0.6888 & 0.2448 & 0.0505 & 0.0119 & 0.0029 & 0.0010 & 0.3984 & 0.4672 \\
        5 & 0.25 & 0.75 & 0.4314 & 0.4500 & 0.0951 & 0.0187 & 0.0039 & 0.0010 & 0.7165 & 0.5711 \\
        5 & 0.75 & 0.25 & 0.1908 & 0.2305 & 0.1540 & 0.1081 & 0.0794 & 0.2373 & 2.3668 & 3.4205 \\
        5 & 0.75 & 0.75 & 0.0330 & 0.2148 & 0.2443 & 0.1665 & 0.1041 & 0.2373 & 2.8060 & 2.4158 \\
        \botrule
    \end{tabular}

    \vspace{1.5cm}

    \caption{CDF $P_{q,\theta}(L_n^{(1)}\leq k)$, for $n=3,5$ and
    $(\theta,q)\in\{0.25,0.75\}^2$.}
    \label{tab:cdf}
    \begin{tabular}{ccccccccc}
        \toprule
        \multirow{2}{*}{$n$} & \multirow{2}{*}{$\theta$} & \multirow{2}{*}{$q$} &
        \multicolumn{6}{c}{$P_{q,\theta}(L_n^{(1)}\leq k)$} \\
        \cmidrule{4-9}
         & & & $k=0$ & $k=1$ & $k=2$ & $k=3$ & $k=4$ & $k=5$ \\
        \midrule
        3 & 0.25 & 0.25 & 0.6921 & 0.9346 & 0.9844 & 1.0000 & --     & --     \\
        3 & 0.25 & 0.75 & 0.5237 & 0.9111 & 0.9844 & 1.0000 & --     & --     \\
        3 & 0.75 & 0.25 & 0.1936 & 0.4287 & 0.5781 & 1.0000 & --     & --     \\
        3 & 0.75 & 0.75 & 0.0632 & 0.3584 & 0.5781 & 1.0000 & --     & --     \\
        \addlinespace
        5 & 0.25 & 0.25 & 0.6888 & 0.9336 & 0.9841 & 0.9961 & 0.9990 & 1.0000 \\
        5 & 0.25 & 0.75 & 0.4314 & 0.8814 & 0.9765 & 0.9952 & 0.9990 & 1.0000 \\
        5 & 0.75 & 0.25 & 0.1908 & 0.4212 & 0.5752 & 0.6833 & 0.7627 & 1.0000 \\
        5 & 0.75 & 0.75 & 0.0330 & 0.2477 & 0.4920 & 0.6586 & 0.7627 & 1.0000 \\
        \botrule
    \end{tabular}
\end{sidewaystable}

The mean $E_{q, \theta}[L_n^{(1)}]$ and variance $\mathrm{Var}_{q, \theta}(L_n^{(1)})$ are obtained using the CDF from the $r$-th moment formula \eqref{eq:moments}.

For instance, for $n=3$, $\theta=q=0.25$, Theorem~\ref{thm:3.1-n} gives $P_{0.25,0.25}(L_3^{(1)}= 0) = 0.6921$, $P_{0.25,0.25}(L_3^{(1)}=1)=0.2424$ and $P_{0.25,0.25}(L_3^{(1)}=2)=0.0498$ (Table~\ref{tab:pmf}), so that $P_{0.25,0.25}(L_3^{(1)}\leq 2)=0.6921+0.2424+0.0498=0.9844$, in agreement with the value obtained directly from Theorem~\ref{thm:3.2-n} and reported in Table~\ref{tab:cdf}.
Moreover, by \eqref{eq:moments} with $r=1$,
\begin{equation*}
    E_{0.25, 0.25}[L_3^{(1)}] = \sum_{k=0}^{2}\big[1-P_{0.25,0.25}(L_3^{(1)}\leq k)\big] =3 - (0.6921 + 0.9346 + 0.9844) = 0.3889,
\end{equation*}
which coincides with the entry $E_{q, \theta}[L_3^{(1)}] = 0.3889$ of Table~\ref{tab:pmf}.

The entries of Tables~\ref{tab:pmf} and~\ref{tab:cdf} illustrate the effect of $\theta$, $q$ and $n$ on the distribution of $L_n^{(1)}$. 
For fixed $q$ and $n$, the expected longest success run length $E_{q, \theta}[L_n^{(1)}]$ increases substantially as $\theta$ increases from $0.25$ to $0.75$, reflecting the larger mass placed on longer runs.
For fixed $\theta$ and $n$, it also increases with $q$, though less markedly than with $\theta$.
For fixed $\theta$ and $q$, the mean increases slightly with $n$, except when $\theta=0.75$, where the increase from $n=3$ to $n=5$ is considerable, as a longer sequence affords substantially more room for extended runs.
Interestingly, $\mathrm{Var}_{q, \theta}(L_n^{(1)})$ does not follow the mean's pattern uniformly in $q$.
For $\theta=0.25$, it increases with $q$ (as does the mean), whereas for $\theta=0.75$ it decreases as $q$ increases from $0.25$ to $0.75$, e.g., from $3.4205$ to $2.4158$ for $n=5$, indicating that the higher success probability concentrates the distribution of $L_n^{(1)}$ around its larger mean rather than spreading it out further.

\begin{figure}[ht]
    \centering
    \includegraphics[width=\textwidth]{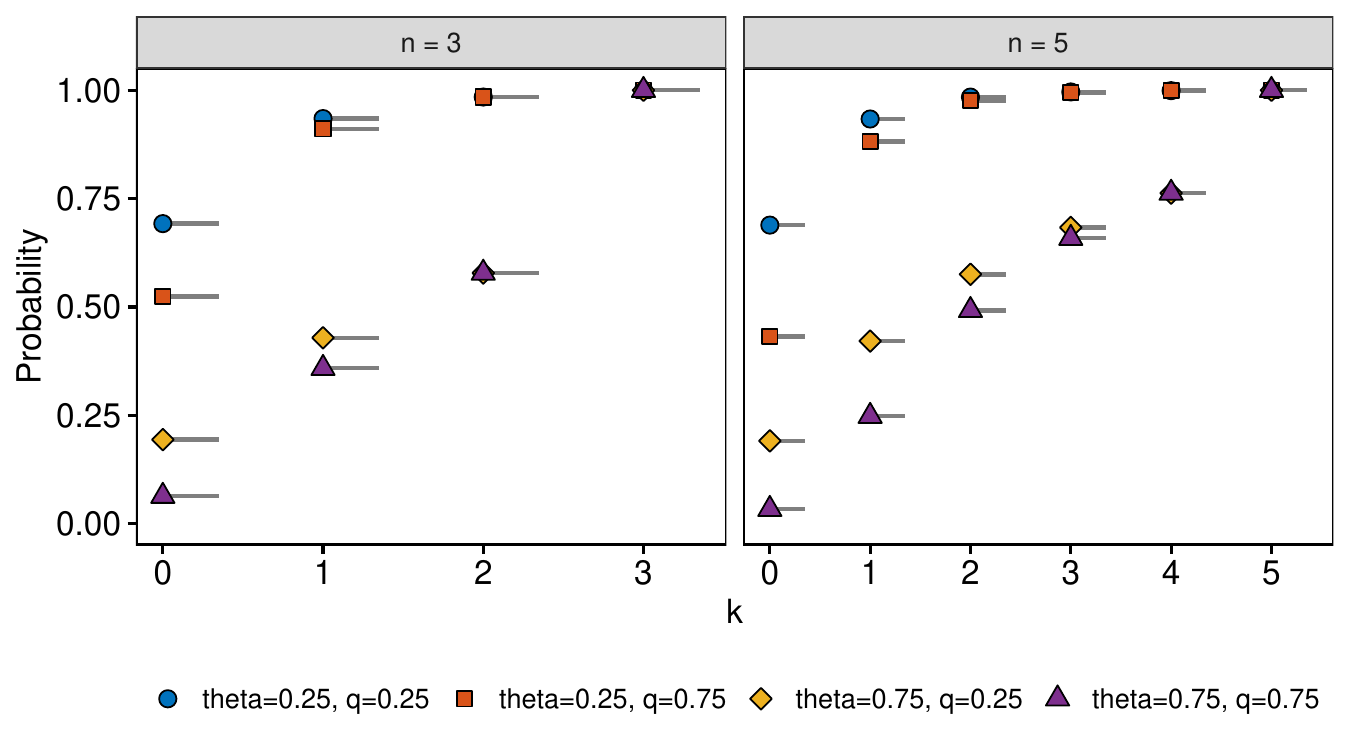}
    \caption{CDF of the longest success run length $L_n^{(1)}$, for $n=3,5$ and $(\theta,q)\in\{0.25,0.75\}^2$.}
    \label{fig:cdf}
\end{figure}

\begin{figure}[ht]
    \centering
    \includegraphics[width=\textwidth]{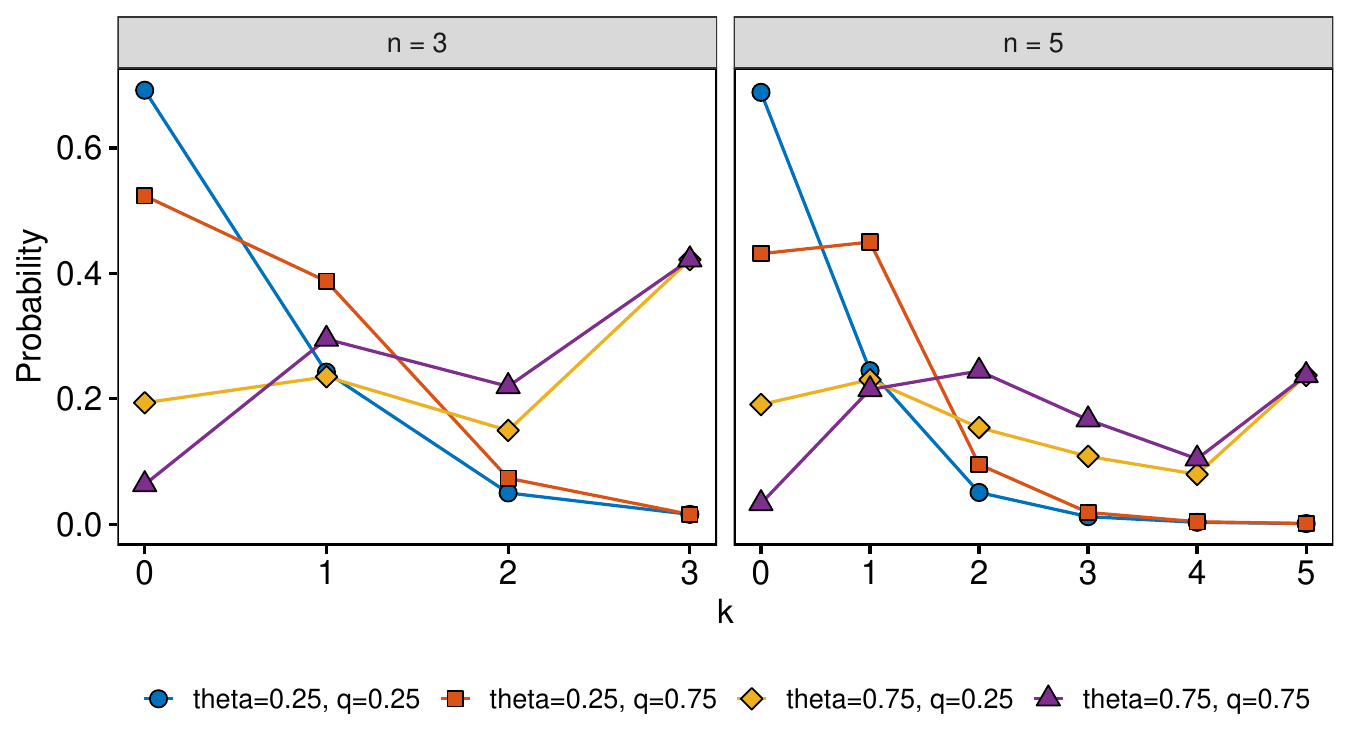}
    \caption{PMF of the longest success run length $L_n^{(1)}$, for $n=3,5$ and $(\theta,q)\in\{0.25,0.75\}^2$. The function is defined only at integer $k$ and the connecting lines serve only as visual guides.}
    \label{fig:pmf}
\end{figure}

Figures~\ref{fig:cdf} and~\ref{fig:pmf} display the CDF and PMF of $L_n^{(1)}$, respectively, for $n=3,5$ and the four combinations of $\theta,q\in\{0.25,0.75\}$ reported in Tables~\ref{tab:pmf} and~\ref{tab:cdf}. 
The plots corroborate the observations above: the CDF curves for $\theta=0.75$ lie visibly below those for $\theta=0.25$ over most of the support of $L_n^{(1)}$, indicating a stochastically larger longest success run length, and the corresponding PMF curves are shifted, and less peaked, toward larger values of $k$.

\section{Simulation study}
\label{sec:sim}

In this section, the distributional results established in Theorems~\ref{thm:3.1-n} and~\ref{thm:3.2-n} are verified by means of a simulation study.
Three aspects are investigated: the accuracy of the maximum likelihood estimator of the parameter vector
$(\theta, q)$, the power of the likelihood ratio test of the hypothesis $H_{0}: q = 1$ against the alternative $H_{1}: q \neq 1$, and the pointwise agreement between the analytical and empirical cumulative distribution functions of $L_{n}^{(1)}$.
The study is conducted over a grid of parameter values $\theta \in \{0.25, 0.75\}$, $q \in \{0.25, 0.75, 1\}$, and sequence lengths $n \in \{10, 20\}$, yielding twelve parameter configurations in total.
For each configuration, $M = 1000$ independent $q$-binary sequences of length $n$ are generated according to Algorithm~\ref{alg:sim} below, and the length of the longest success run $L_{n}^{(1)}$ is recorded from each sequence, yielding a random sample of size $N = 1000$ per replicate dataset.
The configurations with $q \in \{0.25, 0.75\}$ correspond to interior parameter values, while $q = 1$ lies at the boundary of the parameter space $(0, 1]$, where the $q$-binary model reduces to a sequence of independent and identically distributed Bernoulli trials with success probability $\theta$.

\begin{algorithm}
\caption{Simulation of $N$ independent realisations of $L_{n}^{(1)}$ from $q$-binary trials of length $n$}
\label{alg:sim}
    \begin{algorithmic}[1]
    \Require $N \geq 1$, $n \geq 1$, $\theta \in (0,1)$, $q \in (0,1]$
    \For{$i = 1, \ldots, N$}
        \State Set $f \gets 0$, $c \gets 0$, $L \gets 0$
        \For{$j = 1, \ldots, n$}
            \State Generate $U \sim \mathrm{Uniform}(0,1)$
            \If{$U \leq \theta q^{f}$}
                \State Set $X_j \gets 1$, $c \gets c + 1$, $L \gets \max(L, c)$
            \Else
                \State Set $X_j \gets 0$, $c \gets 0$, $f \gets f + 1$
            \EndIf
        \EndFor
        \State Output $\ell_i \gets L$
    \EndFor
    \end{algorithmic}
\end{algorithm}

\subsection{Maximum likelihood estimation}

Let $\ell = (\ell_1, \ell_2, \ldots, \ell_N)$ be an observed random sample of longest success run lengths from a $q$-binary sequence of length $n$.
The log-likelihood function for the parameter vector $(\theta, q)$ is
\begin{equation} \label{eq:loglik}
    \mathcal{L}(\theta, q;\, \ell) = \sum_{i=1}^{N} \log P_{q,\theta}\!\left(L_{n}^{(1)} = \ell_i\right),
\end{equation}
where $P_{q,\theta}(L_{n}^{(1)} = k)$ is given by Theorem~\ref{thm:3.1-n}.
The maximum likelihood estimator $(\hat{\theta}_{ML}, \hat{q}_{ML})$ is the value of $(\theta, q)$ that maximises~\eqref{eq:loglik} over the parameter space $(0,1) \times (0,1]$.
The evaluation of~\eqref{eq:loglik} requires repeated computation of $P_{q,\theta}(L_{n}^{(1)} = k)$, which
depends on the recursive function $U_q^k(r,s,t)$ defined in Lemma~\ref{lemma3.1-n}, posing numerical challenges for maximization.
To reduce the number of likelihood evaluations, the maximization proceeds in two stages: a grid search over a uniform grid on $(0,1)\times(0,1]$ provides an initial estimate, which is subsequently refined by the Limited-memory Broyden--Fletcher--Goldfarb--Shanno algorithm with box constraints \citep{byrd1995limited}.
As no closed form for the gradient of the log-likelihood is available, it is approximated numerically via finite differences with Richardson extrapolation \citep{burden2011numerical}.

For the interior parameter values, $(\theta, q) \in (0,1)\times(0,1)$, the maximum likelihood estimator satisfies the standard regularity conditions and is asymptotically normally distributed \citep{lehmann1998theory},
\begin{equation} \label{eq:asymp_normality}
    \sqrt{N}\left(\begin{pmatrix} \hat{\theta}_{ML} \\ \hat{q}_{ML} \end{pmatrix} - \begin{pmatrix} \theta \\ q \end{pmatrix} \right) \xrightarrow{d} \mathcal{N}_2\!\left(\mathbf{0},\, \mathcal{I}( \theta,q)^{-1}\right), \quad N \to \infty,
\end{equation}
where $\mathcal{I}(\theta, q) = -\mathbb{E}[\nabla^2 \mathcal{L}(\theta, q;\, \ell)]$ is the expected Fisher information matrix.
Asymptotic confidence intervals for $\theta$ and $q$ are therefore constructed only for interior configurations.
They are obtained from the observed Fisher information matrix $\mathcal{J}(\hat{\theta}_{ML}, \hat{q}_{ML}) = -\nabla^2 \mathcal{L}(\hat{\theta}_{ML}, \hat{q}_{ML};\, \ell)$, which is a consistent estimator of $\mathcal{I}(\theta, q)$ \citep{lehmann1998theory} and is approximated numerically via Richardson extrapolation.
Specifically, the asymptotic $95\%$ confidence interval for $\theta$ is $\hat{\theta}_{ML} \pm 1.96\,[\mathcal{J}(\hat{\theta}_{ML}, \hat{q}_{ML})^{-1}]_{11}^{1/2}$, and analogously for $q$, where $[\mathcal{J}(\theta, q)^{-1}]_{jj}$ denotes the $j$-th diagonal entry of the inverse of the observed Fisher information matrix.

At the boundary $q = 1$, the standard regularity conditions for asymptotic normality are not satisfied \citep{lehmann1998theory}, but the MLE remains consistent.

\begin{figure}[htp]
    \centering
    \includegraphics[width=\textwidth]{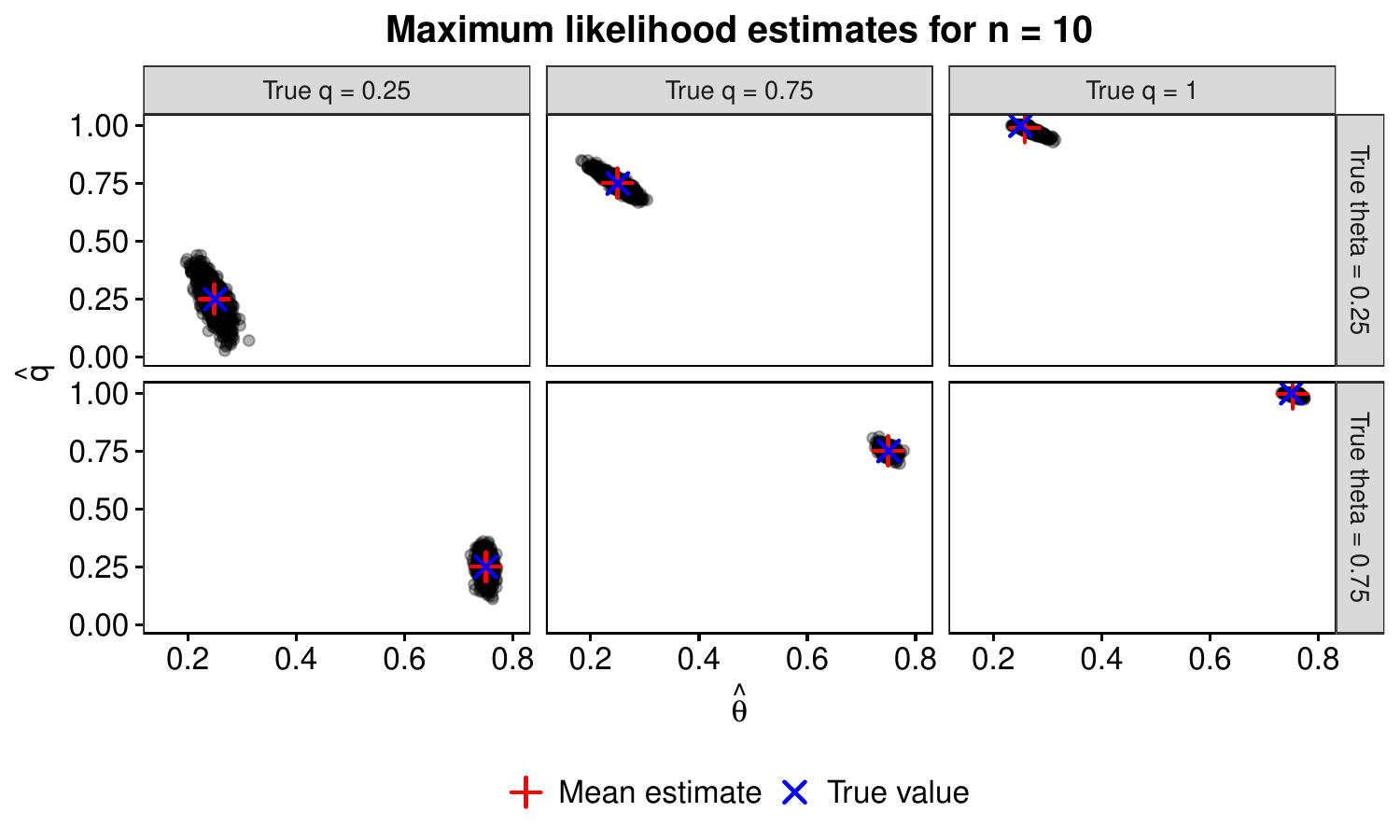}\\[6pt]
    \includegraphics[width=\textwidth]{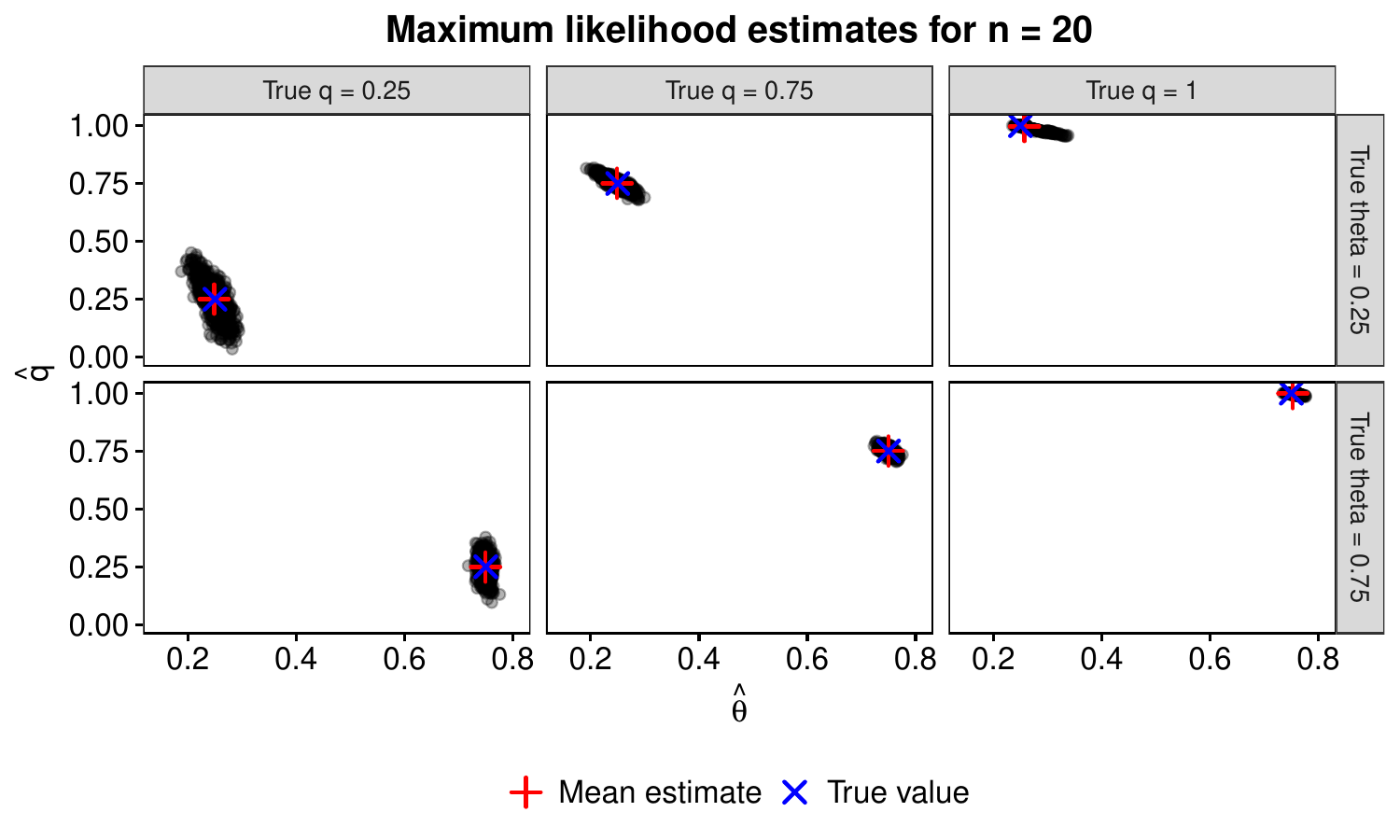}
    \caption{Joint distribution of the maximum likelihood estimates $(\hat{\theta}_{ML}, \hat{q}_{ML})$ across $1000$ replicate datasets, each consisting of $1000$ longest success run lengths recorded from simulated $q$-binary sequences of length $n$. Top: $n=10$; bottom: $n=20$. True parameter values are indicated by the cross markers}
    \label{fig:mle}
\end{figure}

The performance of the estimator is assessed over $M = 1000$ replicate datasets for each of the twelve parameter configurations.
For each configuration, we report the root mean squared error (RMSE)
\begin{equation*}
    \mathrm{RMSE}(\hat{\theta}_{ML}) = \sqrt{\frac{1}{M}\sum_{j=1}^{M}\left(\hat{\theta}_{ML}^{(j)} - \theta\right)^2}.
\end{equation*}
For the configurations in the interior of the parameter space, the empirical coverage of the asymptotic $95\%$ confidence interval is also reported.
The same quantities are computed for $\hat{q}_{ML}$.
The results are presented in Table~\ref{tab:mle} and illustrated in the scatter plots of Figure~\ref{fig:mle}.

\begin{table}[ht]
\caption{Root mean squared error of $\hat{\theta}_{ML}$ and $\hat{q}_{ML}$ and empirical coverage of the asymptotic $95\%$ confidence intervals across $M = 1000$ replicate datasets, each consisting of $N = 1000$ longest success run lengths recorded from simulated $q$-binary sequences of length $n$, for $n \in \{10, 20\}$, $\theta \in \{0.25, 0.75\}$, and $q \in \{0.25, 0.75, 1\}$. Dashes indicate that asymptotic confidence intervals are not constructed at the boundary $q = 1$.}\label{tab:mle}
\centering
    \begin{tabular}{ccccccc}
        \toprule
        $n$ & $\theta$ & $q$ & $\text{RMSE}(\hat{\theta}_{ML})$ & $\text{RMSE}(\hat{q}_{ML})$ & \makecell{Asymp. 95\% CI\\ coverage ($\hat{\theta}_{ML}$)} & \makecell{Asymp. 95\% CI\\ coverage ($\hat{q}_{ML}$)} \\
        \midrule
        10 & 0.25 & 0.25 & 0.0162 & 0.0664 & 0.962 & 0.95 \\
        10 & 0.25 & 0.75 & 0.0193 & 0.0312 & 0.941 & 0.943 \\
        10 & 0.25 & 1 & 0.0159 & 0.0169 & -- & -- \\
        10 & 0.75 & 0.25 & 0.00697 & 0.0409 & 0.956 & 0.947 \\
        10 & 0.75 & 0.75 & 0.00777 & 0.0146 & 0.949 & 0.949 \\
        10 & 0.75 & 1 & 0.0067 & 0.0059 & -- & -- \\
        \addlinespace
        20 & 0.25 & 0.25 & 0.0173 & 0.0701 & 0.947 & 0.937 \\
        20 & 0.25 & 0.75 & 0.0164 & 0.0221 & 0.952 & 0.951 \\
        20 & 0.25 & 1 & 0.0175 & 0.0097 & -- & -- \\
        20 & 0.75 & 0.25 & 0.00684 & 0.0417 & 0.954 & 0.949 \\
        20 & 0.75 & 0.75 & 0.00745 & 0.0136 & 0.955 & 0.934\\
        20 & 0.75 & 1 & 0.0067 & 0.0035 & -- & -- \\
        \botrule
    \end{tabular}
\end{table}

The entries of Table~\ref{tab:mle} show that RMSE is small throughout, including at the boundary $q = 1$, confirming consistency of the MLE across all configurations.
For the interior configurations, $\hat{\theta}_{ML}$ is consistently more precise than $\hat{q}_{ML}$, with the discrepancy most pronounced at $q = 0.25$.
This reflects the fact that strong geometric decay in the success probability is harder to distinguish from moderate decay on the basis of the longest run alone.
The empirical coverage of the asymptotic $95\%$ confidence intervals is close to the nominal level for both parameters across all interior configurations.
As asymptotic normality holds in the limit $N \to \infty$ for interior parameter values, the coverage is not expected to be exact at any finite $N$.
The results confirm that at $N = 1000$ the normal approximation is adequate.
These findings are corroborated by Figure~\ref{fig:mle}, which displays the joint distribution of $(\hat{\theta}_{ML}, \hat{q}_{ML})$ across replicates.

For the interior configurations, the estimates are centred near the true values in all panels, with the spread along the $\hat{q}_{ML}$ axis visibly wider than along the $\hat{\theta}_{ML}$ axis, and the approximately elliptical shape of the point clouds is consistent with the bivariate normal limiting distribution established in~\eqref{eq:asymp_normality}.
For the boundary configurations $q = 1$, the estimates of $\hat{q}_{ML}$ accumulate near the upper boundary of the parameter space, producing the asymmetric, non-elliptical point clouds visible in the rightmost column of Figure~\ref{fig:mle}, as expected when the true parameter lies on the boundary.

\subsection{Likelihood ratio test of \texorpdfstring{$H_0: q = 1$}{H0: q = 1}}

A natural question in practice is whether the observed longest success run lengths are consistent with a standard sequence of independent and identically distributed Bernoulli trials, or whether the geometric variation in the success probability characteristic of the $q$-binary model is supported by the data.
This is formalized as the hypothesis testing problem
\begin{equation*}
    H_0: q = 1 \quad \text{against} \quad H_1: q \neq 1,
\end{equation*}
where under $H_0$ the $q$-binary model reduces to a sequence of IID Bernoulli trials with success probability $\theta$, which remains a free nuisance parameter.
The likelihood ratio test statistic is
\begin{equation*}
    \Lambda = 2\left[\mathcal{L}(\hat{\theta}_{ML}, \hat{q}_{ML};\,\ell) - \mathcal{L}(\hat{\theta}_0, 1;\ \ell)\right],
\end{equation*}
where $\hat{\theta}_0 = \arg\max_{\theta} \mathcal{L}(\theta, 1;\,\ell)$ is the restricted maximum likelihood estimator of $\theta$ under $H_0$, obtained by maximizing the log-likelihood with $q$ fixed at $1$.

Because $q = 1$ lies on the boundary of the parameter space $(0, 1]$, the likelihood ratio $\Lambda$ does not asymptotically follow a chi-squared distribution with degree $1$, denoted as $\chi^2_1$, under $H_0$ as $N \to \infty$ \citep{self1987asymptotic}.
Instead, the null distribution of $\Lambda$ is the equal-weight mixture 
\begin{equation*}
    \Lambda \xrightarrow{d} \frac{1}{2}\chi^2_0 + \frac{1}{2}\chi^2_1, \quad N \to \infty,
\end{equation*}
where $\chi^2_0$ denotes a point mass at zero.
The boundary-corrected $p$-value is accordingly
\begin{equation*}
    p = \frac{1}{2} P\!\left(\chi^2_1 > \Lambda\right) \quad \text{for } \Lambda > 0,
\end{equation*}
and $H_0$ is rejected at the nominal level $\alpha = 0.05$ when $p < 0.05$, which is equivalent to using the $\chi^2_1$ critical value at level $2\alpha = 0.10$.

The power of the test is examined across the eight interior configurations with $q \in \{0.25, 0.75\}$, and the empirical size is verified at the four boundary configurations with $q = 1$.
The results are summarized in Table~\ref{tab:lrt}, which reports the mean and standard deviation of $\Lambda$, the mean and maximum boundary-corrected $p$-value, and the rejection rate across $M = 1000$ replicates.

The empirical rejection rates at $q = 1$ are $0.043$, $0.044$, $0.041$, and $0.048$ across the four null configurations, all close to the nominal level of $0.05$, confirming that the boundary-corrected test is well calibrated.
For the interior configurations, the test achieves high power throughout, with rejection rates between $0.94$ and $1$.
Overall, the results confirm that the likelihood ratio test based on the probability mass function established in Theorem~\ref{thm:3.1-n} is correctly sized under $H_0$ and achieves high power against $H_1: q \neq 1$, as expected.

\begin{table}[ht]
\caption{Summary of the likelihood ratio test of $H_0: q = 1$ against $H_1: q \neq 1$ across $M = 1000$ replicate datasets, each consisting of $N = 1000$ longest success run lengths recorded from simulated $q$-binary sequences of length $n$, for $n \in \{10, 20\}$, $\theta \in \{0.25, 0.75\}$, and $q \in \{0.25, 0.75, 1\}$. $p$-values are boundary-corrected following \citet{self1987asymptotic}. Size rows correspond to $q = 1$; power rows correspond to $q < 1$. Rejection is at the nominal level $\alpha = 0.05$}\label{tab:lrt}
\centering
    \begin{tabular}{cccccccc}
        \toprule
        $n$ & $\theta$ & $q$ & Mean $\Lambda$ & SD$(\Lambda)$ & Mean $p$-value & Max $p$-value & \makecell{Rejection \\ rate} \\
        \midrule
        10 & 0.25 & 0.25 & 165.9 & 34.0  & $8.93 \cdot 10^{-22}$ & $5.49 \cdot 10^{-19}$ & 1 \\
        10 & 0.25 & 0.75 & 50.4 & 16.3  & $9.39 \cdot 10^{-7}$  & $3.91 \cdot 10^{-4}$  & 1 \\
        10 & 0.25 & 1 & 0.4 & 1  & 0.638 & 1 & 0.043 \\
        10 & 0.75 & 0.25 & 2208.4 & 159.3 & 0 & 0 & 1 \\
        10 & 0.75 & 0.75 & 571.7 & 102.8 & $3.16 \cdot 10^{-60}$ & $3.01 \cdot 10^{-57}$ & 1 \\
        10 & 0.75 & 1 & 0.5 & 1.1 & 0.628 & 1 & 0.044 \\
        \addlinespace
        20 & 0.25 & 0.25 & 255.2 & 105.9 & 0.011 & 1 & 0.989 \\
        20 & 0.25 & 0.75 & 136.8 & 84.6 & 0.057 & 1 & 0.942 \\
        20 & 0.25 & 1 & 0.3 & 1.1 & 0.698 & 1 & 0.041 \\
        20 & 0.75 & 0.25 & 3552.5 & 253.0 & 0 & 0 & 1 \\
        20 & 0.75 & 0.75 & 1112.3 & 165.6 & $7.72 \cdot 10^{-145}$& $7.68 \cdot 10^{-142}$ & 1 \\
        20 & 0.75 & 1 & 0.5 & 1.1 & 0.640 & 1 & 0.048 \\
        \botrule
    \end{tabular}
\end{table}

\subsection{Comparison of analytical and empirical cumulative
distribution functions}

The cumulative distribution function established in Theorem~\ref{thm:3.2-n} is assessed by comparing it pointwise against the empirical cumulative distribution function of $L_{n}^{(1)}$ constructed from the pooled sample of $N \cdot M = 10^6$ observations for each of the twelve parameter configurations.
Let $F_0(k; q, \theta) = P_{q,\theta}(L_{n}^{(1)} \leq k)$ denote the analytical CDF given by Theorem~\ref{thm:3.2-n} and $F_N(k) = (NM)^{-1} \sum_{i=1}^{NM} \mathbf{1}(L_{n,i}^{(1)} \leq k)$ the corresponding empirical CDF.
Since the support of $L_{n}^{(1)}$ is the finite set $\{0, 1, \ldots, n\}$, pointwise convergence of $F_N(k)$ to $F_0(k; q, \theta)$ at each $k$ is equivalent to uniform convergence, and the Kolmogorov-Smirnov statistic
\begin{equation*}
    D = \sup_{k \in \{0,\ldots,n\}} |F_N(k) - F_0(k; q, \theta)|
\end{equation*}
therefore reduces to a maximum over finitely many points.
As the Kolmogorov-Smirnov statistic is not distribution-free for discrete distributions, $D$ is reported as a descriptive measure of the maximum pointwise discrepancy only, without performing the Kolmogorov-Smirnov test of equality of CDFs.

Instead, formal inference is based on exact two-sided binomial tests.
Let $G(k) = P(L_{n}^{(1)} \leq k)$ denote the true cumulative distribution function of $L_{n}^{(1)}$.
The pointwise hypothesis
\begin{equation}
    \label{eq:cdf_hyp}
    H_0^{(k)}: G(k) = F_0(k; q, \theta) \quad \text{against} \quad H_1^{(k)}: G(k) \neq F_0(k; q, \theta)
\end{equation}
is examined for each interior support point $k \in \{0, 1, \ldots, n-1\}$.

Under $H_0^{(k)}$, each observation $L_{n,i}^{(1)}$ for $i = 1, \dots, NM$ satisfies $P(L_{n,i}^{(1)} \leq k) = F_0(k; q, \theta)$, so the indicators $\mathbf{1}(L_{n,i}^{(1)} \leq k)$ are independent Bernoulli random variables with success probability $F_0(k; q, \theta)$.
Consequently, the count
\begin{equation*}
    S_k = \sum_{i=1}^{NM} \mathbf{1}\!\left(L_{n,i}^{(1)} \leq k\right)
\end{equation*}
follows a binomial distribution with parameters $NM$ and $F_0(k; q, \theta)$ under $H_0^{(k)}$, from which an exact two-sided $p$-value is obtained.
The point $k = n$ is excluded as $P_{q,\theta}(L_{n}^{(1)} \leq n) = 1$ holds trivially for any distribution on $\{0, \ldots, n\}$.
Because $n$ hypotheses of the form~\eqref{eq:cdf_hyp} are tested simultaneously on the same pooled sample, the tests are not mutually independent.
The Bonferroni correction controls the family-wise error rate at level $\alpha = 0.05$ without requiring independence of the individual tests \citep{lehmann2022multiple}, and the family-wise rejection threshold for each individual test is accordingly $\alpha / n$ with $\alpha = 0.05$.

\begin{figure}[htp]
    \centering
    \includegraphics[width=\textwidth]{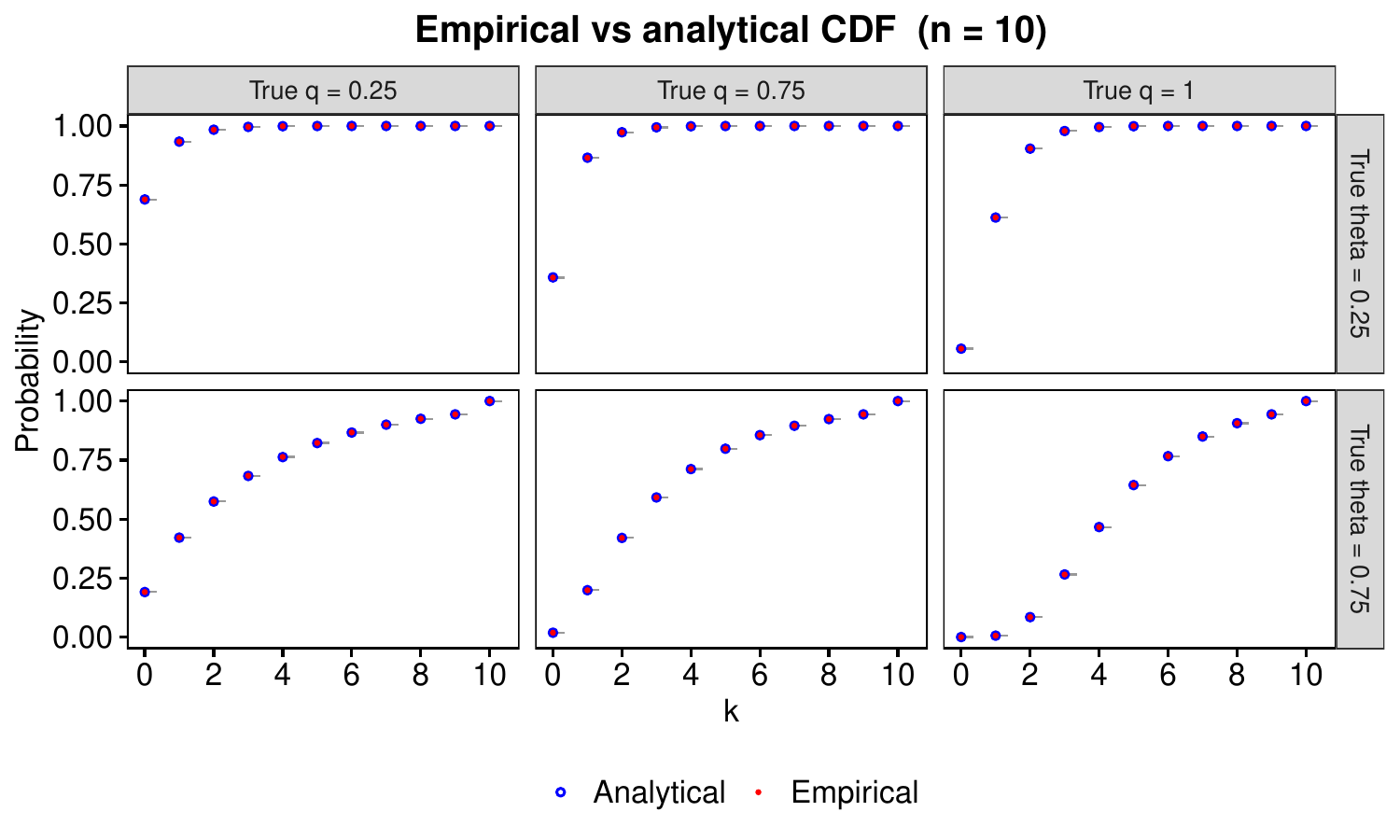}\\[6pt]
    \includegraphics[width=\textwidth]{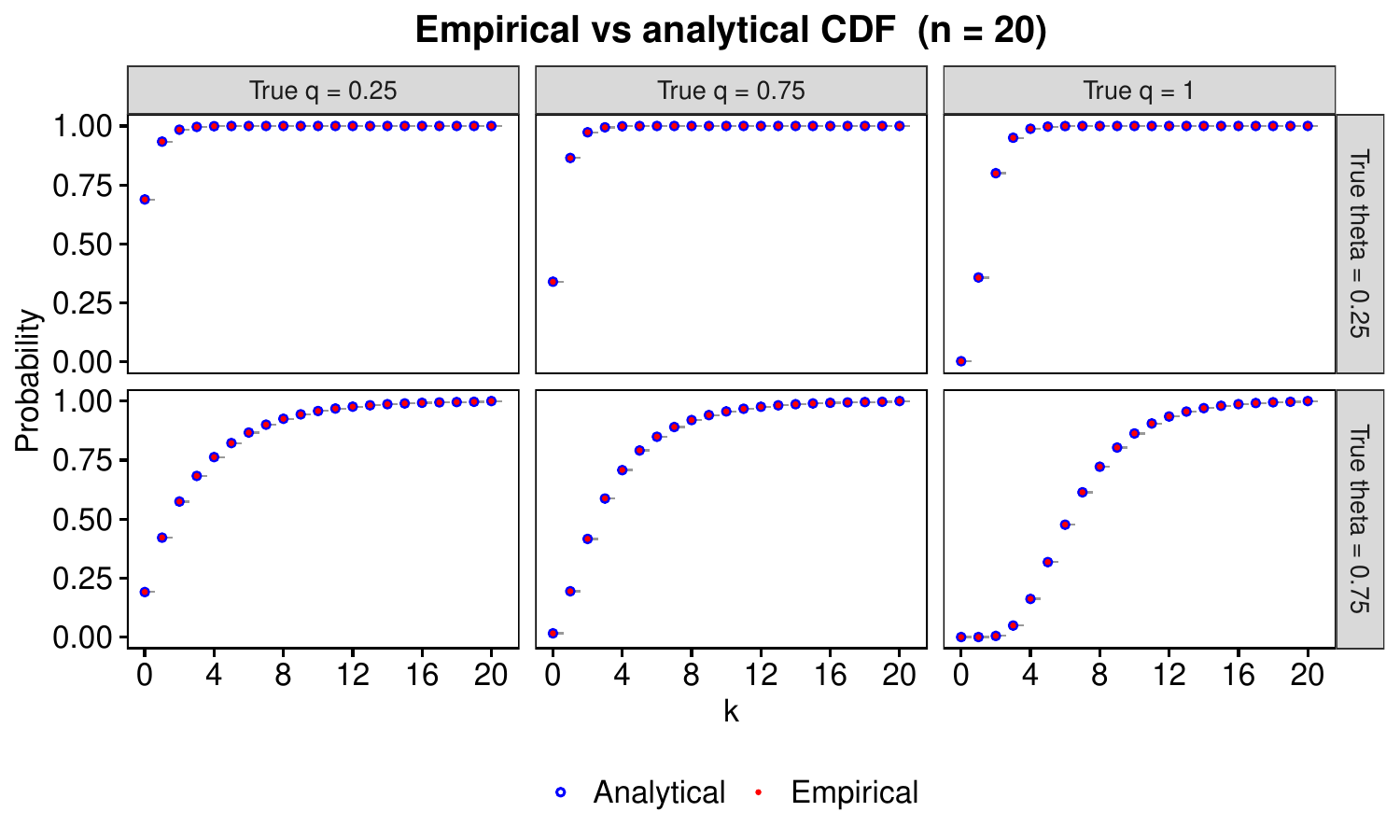}
    \caption{Analytical (blue, hollow) and empirical (red, solid) cumulative distribution functions of $L_{n}^{(1)}$ for $n = 10$ (top) and $n = 20$ (bottom), constructed from $NM = 10^6$ observations per configuration. The two functions are visually indistinguishable in all panels}
    \label{fig:cdf_test}
\end{figure}

The results are presented in Table~\ref{tab:cdf_test} and the analytical and empirical CDFs are displayed in Figure~\ref{fig:cdf_test}.
The null hypothesis $H_0^{(k)}$ is not rejected at any support point in any of the twelve configurations. The values of $D$ are of order $10^{-4}$ to $10^{-3}$ throughout, reflecting negligible pointwise discrepancy between the analytical and empirical CDFs at the pooled sample size of $10^6$.
These results confirm the correctness of the cumulative distribution function established in Theorem~\ref{thm:3.2-n} across all parameter configurations considered.

\begin{table}[ht]
\caption{Bonferroni-corrected exact binomial test results and maximum pointwise discrepancy $D = \sup_k |F_N(k) - F_0(k; q, \theta)|$ between the empirical CDF constructed from $NM = 10^6$ observations and the analytical CDF of Theorem~\ref{thm:3.2-n}, for $n \in \{10, 20\}$, $\theta \in \{0.25, 0.75\}$, and $q \in \{0.25, 0.75, 1\}$. The null hypothesis is not rejected in any configuration at $\alpha = 0.05$}\label{tab:cdf_test}
\centering
    \begin{tabular}{cccccc}
        \toprule
        $n$ & $\theta$ & $q$ & Bonferroni threshold $(\alpha/n)$ &
        Min $p$-value & $D$ \\
        \midrule
        10 & 0.25 & 0.25 & 0.0050 & 0.0083 & $2.90 \cdot 10^{-4}$ \\
        10 & 0.25 & 0.75 & 0.0050 & 0.0373 & $9.99 \cdot 10^{-4}$ \\
        10 & 0.25 & 1    & 0.0050 & 0.2071 & $3.07 \cdot 10^{-4}$ \\
        10 & 0.75 & 0.25 & 0.0050 & 0.0472 & $8.45 \cdot 10^{-4}$ \\
        10 & 0.75 & 0.75 & 0.0050 & 0.4198 & $2.84 \cdot 10^{-4}$ \\
        10 & 0.75 & 1    & 0.0050 & 0.2765 & $3.39 \cdot 10^{-4}$ \\
        \addlinespace
        20 & 0.25 & 0.25 & 0.0025 & 0.0560 & $2.23 \cdot 10^{-4}$ \\
        20 & 0.25 & 0.75 & 0.0025 & 0.1683 & $4.18 \cdot 10^{-4}$ \\
        20 & 0.25 & 1    & 0.0025 & 0.0042 & $1.15 \cdot 10^{-3}$ \\
        20 & 0.75 & 0.25 & 0.0025 & 0.0173 & $7.21 \cdot 10^{-4}$ \\
        20 & 0.75 & 0.75 & 0.0025 & 0.0818 & $7.92 \cdot 10^{-4}$ \\
        20 & 0.75 & 1    & 0.0025 & 0.1112 & $5.87 \cdot 10^{-4}$ \\
        \botrule
    \end{tabular}
\end{table}

\section{Discussion}\label{sec:discussion}

In the present work, we introduce and study the distribution of the length of the longest success run for $q-$binary trials. We have obtained the exact expression for the PMF of the random variable $L_{n}^{(1)}$ which follows a $q-L^{(1)}_{n}(\theta)$ distribution. More specifically, we derive the exact expression for the PMF $P_{q,\theta}\left(L_{n}^{(1)}= k\right)$, representing the distribution of the length $L_{n}^{(1)}$ of the longest success run and the cumulative distribution function $P_{q,\theta}\left(L_{n}^{(1)}\leq k\right)$ for $L_{n}^{(1)}$ through enumerative combinatorics. These theoretical results are further validated by a simulation study, which confirms the accuracy of maximum likelihood estimation of $(\theta, q)$, the correct size and high power of a boundary-corrected likelihood ratio test of $H_0: q=1$, and close agreement between the derived analytical cumulative distribution function of $L_{n}^{(1)}$ and its empirical counterpart.

Another random variable that is closely associated with $L_{n}^{(1)}$ is the $T_{k}$, which represents the waiting time until the first occurrence of a sequence of $k$ consecutive successes. In the context of a $q-$geometric model, the distribution of the $T_{k}$ is referred to as $q-$geometric distribution of order $k$, can be denoted as $q-T_{k}(\theta)$. Given that the event $\left\{L_{n}^{(1)} < k\right\}$ is equivalent to $\Big\{T_{k} > n\Big\}$, we can directly establish the following relationship
\noindent
\begin{equation}
\label{dual}
P_{q,\theta}\left(L_{n}^{(1)}< k\right)=P_{q,\theta}\Big(T_{k} > n\Big).
\end{equation}
\noindent
This relationship provides straightforward methods for determining the exact distribution (probability mass function) of $T_{k}$ by utilizing the formulae derived in this study for the longest run $L_{n}^{(1)}$. The probability of the event $\{T_{k} = n\}$ can be expressed as follows.
\noindent
\begin{equation*}
    P_{q,\theta}\Big(T_{k} = n\Big)=P_{q,\theta}\Big(T_{k} > n-1\Big)-P_{q,\theta}\Big(T_{k}> n\Big).
\end{equation*}
\noindent
Using dual relationship \eqref{dual}, we can rewrite as follows.
\noindent
\begin{equation*}
    P_{q,\theta}\Big(T_{k} = n\Big)=P_{q,\theta}\Big(L_{n-1}^{(1)}< k\Big)-P_{q,\theta}\Big(L_{n}^{(1)}< k\Big).
\end{equation*}
\noindent
Using simple mathematical algebra arguments, we can rewrite as follows.
\noindent
\begin{equation}
\label{eq5.2}
    P_{q,\theta}\Big(T_{k} = n\Big)
    =P_{q,\theta}\Big(L_{n-1}^{(1)}\leq k-1\Big)-P_{q,\theta}\Big(L_{n}^{(1)}\leq k-1\Big).
\end{equation}
\noindent
Using Formula \eqref{eq5.2} and Theorem \ref{thm:3.2-n}, we can compute the PMF for a $q-$geometric distribution of order $k$. Alternatively, we can compute the PMF values for the $q-$geometric distribution of order $k$ using the formulae derived by \citet{yalcin2014q} in a manner similar to that used in Theorem \ref{thm:3.2-n}.\\
As a possible direction for future research, one could consider deriving the recursive expressions for the PMF of the RV $T_{k}$. The preliminary results from the investigation of this RV are promising. To give a specific example, the PMF of $T_k$, $f_q(x;\theta)=P_{q,\theta}\left(T_{k} = x\right),$ $x=k,k+1,\ldots$ satisfies for $x\geq k+1$ the recursive scheme
\noindent
\begin{equation*}
    f_q(x;\theta)
    =\sum_{i=1}^{k}\theta^{i-1}(1-\theta)f_q(x-i;\theta q),\ \text{and}
\end{equation*}

\begin{equation*}
    f_q(k;\theta)
    =\theta^{k},\ \text{for}\ x=k,\ \text{and}
\end{equation*}

\begin{equation*}
    f_q(x;\theta)
    =0,\ \text{for}\ x<k.
\end{equation*}
\noindent


\backmatter

\bmhead{Acknowledgements}

The authors would like to thank Dr.\ habil.\ Tommy Rene Jensen, whose comments led to significant improvements in this manuscript. The authors would like to thank Professor Frosso S. Makri for drawing our attention to the connection between the longest success run and the number of success runs exceeding a fixed threshold. The first author is grateful to Prof. Dr. Mihaela Zavolan and her research group at the Biozentrum, University of Basel, for hosting a research stay that significantly influenced the computational work presented in this paper.

\section*{Declarations}

\begin{itemize}
\item \textbf{Funding} This work was supported by the Ministry of Education of the Republic of Korea and the National Research Foundation of Korea.(NRF-2025S1A5B5A16008020)
\item \textbf{Conflict of interest} The authors have no competing interests to declare that are relevant to the content of this article.
\item \textbf{Ethics approval} Not applicable.
\item \textbf{Consent to participate} Not applicable.
\item \textbf{Consent for publication} Not applicable.
\item \textbf{Data availability} No external datasets were used. The numerical results are based on simulated $q$-binary sequences generated as described in the simulation study.
\item \textbf{Materials availability} Not applicable.
\item \textbf{Code availability} Computational code used to evaluate the derived formulae and to run the simulation study will be made available in a public repository upon acceptance of the manuscript.
\item \textbf{Author contribution} Jungtaek Oh: conceptualization, methodology, theoretical derivation, writing (original draft), review and editing. Nurzhan Sapargali: simulation study, numerical experiments, review and editing. Both authors read and approved the final manuscript.
\end{itemize}

\bibliography{biblio}

\end{document}